\documentclass[12pt,reqno]{article}

\usepackage[usenames]{color}
\usepackage[colorlinks=true,
linkcolor=webgreen, filecolor=webbrown,
citecolor=webgreen]{hyperref}

\definecolor{webgreen}{rgb}{0,.5,0}
\definecolor{webbrown}{rgb}{.6,0,0}

\usepackage{amssymb}
\usepackage{graphicx}
\usepackage{amscd}
\usepackage{lscape}
\usepackage{tikz}
\usepackage{tikz-cd}
\usepackage{pgfplots}

\usetikzlibrary{matrix}
\usetikzlibrary{fit,shapes}
\usetikzlibrary{positioning, calc}
\tikzset{circle node/.style = {circle,inner sep=1pt,draw, fill=white},
        X node/.style = {fill=white, inner sep=1pt},
        dot node/.style = {circle, draw, inner sep=5pt}
        }
\usepackage{tkz-fct}

\usepackage{amsthm}
\newtheorem{theorem}{Theorem}

\theoremstyle{definition}

\newtheorem{example}[theorem]{Example}

\usepackage{float}

\usepackage{graphics,amsmath}
\usepackage{amsfonts}
\usepackage{latexsym}
\usepackage{epsf}

\begin{document}

\begin{center}
\vskip 1cm{\LARGE\bf Power series, the Riordan group and Hopf algebras} \vskip 1cm \large
Paul Barry\\
School of Science\\
Waterford Institute of Technology\\
Ireland\\
\href{mailto:pbarry@wit.ie}{\tt pbarry@wit.ie}
\end{center}
\vskip .2 in

\begin{abstract} The Riordan group, along with its constituent elements, Riordan arrays,  has been a tool for combinatorial exploration since its inception in 1991. More recently, this group has made an appearance in the area of mathematical physics, where it can be used as a toy model in the theory of the renormalization of scalar fields. In this context, its Hopf algebra nature is of importance. In this note, we explain these notions. Power series play a fundamental role in this discussion.   \end{abstract}

\section{The Riordan group}

In this note, we use the notation
$$\mathbb{N} = \{1,2,3,\ldots\}$$
for the set of natural or counting numbers,
$$\mathbb{N}_0 = \{0, 1,2,3,\ldots\}$$ for the set of non-negative integers, and
$$\mathbb{Z} = \{\ldots, -3, -2, -1, 0, 1, 2, 3, \ldots\}$$ for the
set of integers.
The notation $\mathbb{K}$ will be used to denote a field of characteristic $0$, and $\mathbb{R}$ and $\mathbb{C}$ will denote the fields of real and complex numbers, respectively.

Many of the sets that we shall study will be algebras, so we recall the definition of an algebra. Thus we say that a $\mathbb{K}$-algebra is a ring $A$ with unit $1_A$ together with a ring homomorphism $\lambda_A: \mathbb{K} \rightarrow A$ which satisfies $\lambda_A(r)a=a\lambda_A(r)$ for $a \in A, r \in \mathbb{K}$. Then $A$ is a vector space over $\mathbb{K}$ with scalar multiplication given by
$$ra = \lambda_A(r)a = a \lambda_A(r)$$ for $r \in \mathbb{K}, a \in A$, and a product $ m_A: A \times A \rightarrow A$, with $m_A(a,b)=ab$.

A \emph{(formal) power series} over the field $\mathbb{K}$ is a formal expression of the form
$$g(x)=g_0 + g_1 x + g_2 x^2 + g_3 x^3 + \cdots = \sum_{n=0}^{\infty} g_n x^n, \quad g_n \in \mathbb{K}.$$
Such formal power sequences are in a one-to-one correspondence with sequences
$$ g_0, g_1, g_2, g_3, \ldots$$ which can be regarded as maps from $\mathbb{N}_0$ to $\mathbb{K}$. The power series $\sum_{n=0}^{\infty} g_n x^n$ is called the \emph{(ordinary) generating function} of the sequence
$$ g_0, g_1, g_2, g_3, \ldots$$
Using the notations
$$\mathcal{F}=\mathcal{F}(\mathbb{K})=\mathbb{K}[[x]]$$ for the set of formal power series over $\mathbb{K}$, we see that $$ \mathcal{F} \cong \mathbb{K}^{\mathbb{N}_0}$$ by the correspondence
$$ \sum_{n=0}^{\infty} g_n x^n \mapsto (g_0, g_1, g_2, g_3, \ldots).$$

In the above, it is seen that $x$ is a ``dummy'' or ``synthetic'' variable, in that $\sum_{n=0} g_n t^n$ also represents the element $(g_0, g_1, g_2, \ldots)$.

We shall also use the notation $\mathbb{K}[x]$ to denote the algebra of polynomials in the indeterminate $x$ over $\mathbb{K}$. Its elements are thus formal sums of the form $P_n(x)=\sum_{k=0}^n p_k x^k$. The degree of such a polynomial is the highest value of $k$ for which $p_k$ is non-zero. The product of a polynomial of degree $n$ times a polynomial of degree $m$ is a polynomial of degree $n+m$.

The \emph{Riordan group} $\mathcal{R}$ was first defined \cite{SGWW} in 1991 by Shapiro, Getu, Woan and Woodson. As a group of matrices, its elements are invertible lower-triangular matrices with elements
$$d_{n,k}= [x^n] g(x)\phi(x)^k,$$ for suitable power series $g(x)$ and $\phi(x)$ defined over an appropriate ring or field.

Here, $[x^n]$ is the functional

$$ [x^n] : \mathcal{F}=\mathbb{K}[[x]] \longrightarrow \mathbb{K}$$
 $$ f(x)=\sum_{n=0}^{\infty} f_n x^n \mapsto f_n = \frac{1}{n!}\frac{d^n}{dx^n}f|_{x=0}$$ that extracts from the power series $f(x)$ the coefficient of $x^n$. (See the Appendix for more about this functional \cite{MC}).
We note that  $x$ here is a ``dummy'' or ``synthetic'' variable. Hence we have
$$f_n = [x^n] f(x) = [y^n] f(y).$$

\begin{example} The element $\left(\frac{1}{(1-x)^2}, x(1+x)\right) \in \mathcal{R}$ has a matrix representation given by $$\left(
\begin{array}{ccccccc}
 1 & 0 & 0 & 0 & 0 & 0 & \cdots \\
 2 & 1 & 0 & 0 & 0 & 0 & \cdots \\
 3 & 3 & 1 & 0 & 0 & 0 & \cdots \\
 4 & 5 & 4 & 1 & 0 & 0 & \cdots \\
 5 & 7 & 8 & 5 & 1 & 0 & \cdots \\
 6 & 9 & 12 & 12 & 6 & 1 & \cdots \\
 \vdots &  &  &  &  & \vdots & \ddots \\
\end{array}
\right)$$
\end{example}

Many examples of Riordan arrays are to be found in the On-Line Encyclopedia of Integer Sequences \cite{SL1, SL2}. The nature of the matrix representation of the elements of the Riordan group is explored in \cite{EulerConst}.

In order to ensure that these matrices are lower-triangular and invertible, we stipulate that
$$\phi(x)= \phi_1 x + \phi_2 x^2 + \phi_3 x^3+ \cdots,$$ with $\phi_1 \ne 0$. In other words, we have $\phi_0 = 0, \phi_1 \ne 0$.
The \emph{monic} Riordan group $\mathcal{R}^{(1)}$, whose matrices have all $1$'s on the diagonal, can then be prescribed by taking $$g(x)=1 + g_1 x+ g_2 x^2 + \cdots,$$ along with $\phi_1=1$. This is a subgroup of $\mathcal{R}$, being evidently closed under matrix multiplication.

Taken over a field $\mathbb{K}$ of characteristic $0$, the matrices corresponding to $\mathcal{R}^{(1)}$  form a closed subgroup of the prounipotent group \cite{Lubotsky} $T_{\infty}$ of lower-triangular matrices all of whose diagonal elements are $1$. This is a Lie group, thus the monic Riordan group is a Lie subgroup of this group. The Lie algebra $\mathfrak{t}_{\infty}$ of $T_{\infty}$ is composed of the lower-triangular nilpotent matrices with $0$ on the diagonal. The corresponding Lie algebra $\mathfrak{r}^{(1)}$ is composed of those nilpotent matrices of the form $u_{\phi}+d_{\psi}$, for $\phi, \psi \in \mathcal{F}_1$, where \cite{Bacher}
$$u_{\phi}=\left(
\begin{array}{ccccccc}
 \phi_0 & 0 & 0 & 0 & 0 & 0 & \cdots \\
 \phi_1 & \phi_0 & 0 & 0 & 0 & 0 & \cdots \\
 \phi_2 & \phi_1 & \phi_0 & 0 & 0 & 0 & \cdots \\
 \phi_3 & \phi_2 & \phi_1 & \phi_0 & 0 & 0 & \cdots \\
 \phi_4 & \phi_3 & \phi_2 & \phi_1 & \phi_0 & 0 & \cdots \\
 \phi_5 & \phi_4 & \phi_3 & \phi_2 & \phi_1 & \phi_0 & \cdots \\
 \vdots & \vdots & \vdots & \vdots & \vdots & \vdots & \ddots \\
\end{array}
\right),$$ and
$$u_{\psi}=\left(
\begin{array}{ccccccc}
 0 & 0 & 0 & 0 & 0 & 0 & \cdots \\
 0 & \psi_0 & 0 & 0 & 0 & 0 & \cdots \\
 0 & \psi_1 & 2\psi_0 & 0 & 0 & 0 & \cdots \\
 0 & \psi_2 & 2\psi_1 & 3\psi_0 & 0 & 0 & \cdots \\
 0 & \psi_3 & 2\psi_2 & 3\psi_1 & 4\psi_0 & 0 & \cdots \\
 0 & \psi_4 & 2\psi_3 & 3\psi_2 & 4\psi_1 & 5\psi_0 & \cdots \\
\vdots & \vdots & \vdots & \vdots & \vdots & \vdots & \ddots \\
\end{array}
\right).$$
Note that $\phi_0=\psi_0=0$.

For combinatorial purposes, it is often sufficient to consider ``monic'' matrices with integer entries \cite{Barry}. However we shall continue here to consider entities defined over the field $\mathbb{K}$.

\begin{example} The element $\left(\frac{1}{1-x}, \frac{x}{1-x}\right) \in \mathcal{R}$ has general term
\begin{eqnarray*} b_{n,k}&=&[x^n] \frac{1}{1-x} \left(\frac{x}{1-x}\right)^k \\
&=& [x^n] \frac{x^k}{(1-x)^{k+1}}\\
&=& [x^{n-k}] (1-x)^{-(k+1)}\\
&=& [x^{n-k}] \sum_{i=0}^{\infty} \binom{-(k+1)}{i} (-x)^i\\
&=& [x^{n-k}] \sum_{i=0}^{\infty} \binom{k+1+i-1}{i} (-1)^i (-x)^i\\
&=& [x^{n-k}] \sum_{i=0}^{\infty} \binom{k+i}{i} x^i\\
&=& \binom{k+n-k}{n-k} \\
&=& \binom{n}{k}.\end{eqnarray*}

Thus the Riordan array element $\left(\frac{1}{1-x}, \frac{x}{1-x}\right)$ corresponds to the binomial matrix that begins
$$\left(
\begin{array}{ccccccc}
 1 & 0 & 0 & 0 & 0 & 0 & \cdots \\
 1 & 1 & 0 & 0 & 0 & 0 & \cdots  \\
 1 & 2 & 1 & 0 & 0 & 0 & \cdots \\
 1 & 3 & 3 & 1 & 0 & 0 & \cdots \\
 1 & 4 & 6 & 4 & 1 & 0 & \cdots \\
 1 & 5 & 10 & 10 & 5 & 1 & \cdots \\
 \vdots & \vdots & \vdots & \vdots & \vdots & \vdots & \ddots \\
\end{array}
\right).$$

The exponential map
$$\exp : \mathfrak{r} \longrightarrow \mathcal{R} $$ then gives us

$$ \exp\left\{\left(
\begin{array}{ccccccc}
 0 & 0 & 0 & 0 & 0 & 0 & \cdots \\
 1 & 0 & 0 & 0 & 0 & 0 & \cdots \\
 0 & 2 & 0 & 0 & 0 & 0 & \cdots \\
 0 & 0 & 3 & 0 & 0 & 0 & \cdots \\
 0 & 0 & 0 & 4 & 0 & 0 & \cdots \\
 0 & 0 & 0 & 0 & 5 & 0 & \cdots \\
 \vdots & \vdots & \vdots & \vdots & \vdots & \vdots & \ddots \\
\end{array}
\right)\right\} = \left(
\begin{array}{ccccccc}
 1 & 0 & 0 & 0 & 0 & 0 & \cdots \\
 1 & 1 & 0 & 0 & 0 & 0 & \cdots \\
 1 & 2 & 1 & 0 & 0 & 0 & \cdots \\
 1 & 3 & 3 & 1 & 0 & 0 & \cdots \\
 1 & 4 & 6 & 4 & 1 & 0 & \cdots \\
 1 & 5 & 10 & 10 & 5 & 1 & \cdots \\
 \vdots & \vdots & \vdots & \vdots & \vdots & \vdots & \ddots \\
\end{array}
\right).$$

\end{example}

As it happens, \cite{Luzon} the Riordan group $\mathcal{R}$ can be shown to be the inverse limit of the groups $\mathcal{R}_n$, where the group $\mathcal{R}_n$ is obtained by taking the $n \times n$ truncations of the elements of $\mathcal{R}$.

Each element of $\mathcal{R}$ is thus determined by a pair of power series $g(x), \phi(x)$ where
$$ g(x)=1+g_1 x + g_2 x^2+ g_3 x^3 + \cdots$$ and
$$\phi(x) = x +\phi_2 x^2+ \phi_3 x^3+ \cdots$$

We let $\mathcal{F}=\mathbb{K}[[x]]$, and we let
$$\mathcal{F}_d=\{ f(x) \in \mathcal{F} | f(x)=f_d x^d+f_{d+1}x^{d+1} + \cdots, f_d \ne 0\}.$$
The elements of $\mathcal{F}_d$ are called powers series of \emph{order} $d$.
The set $\mathcal{F}$ is an algebra for the convolution product
$$ f(x)g(x)=f(x).g(x) = \sum_{n=0}^{\infty} (\sum_{k=0}^n f_k g_{n-k}) x^n,$$ and the
scalar product $$\mathbb{K} \times \mathcal{F} \rightarrow \mathcal{F},$$
$$ (r, g) \mapsto rg$$
The identity is the element
$$1 = 1.x^0 + 0 x^1 + 0 x^2 + \cdots,$$ corresponding to the element
$$(1, 0,0,0,\ldots) \in \mathbb{K}^{\mathbb{N}_0}.$$
This convolution product defines a multiplication
$$m_c: \mathcal{F} \times \mathcal{F} \longrightarrow \mathcal{F},$$
$$     (g, f) \mapsto gf,$$ which is commutative.

The set $\mathcal{F}_0=\{ g(x) \in \mathcal{F} | g(x)=g_0 +g_1x + \cdots, g_0 \ne 0\}$ is the group of \emph{invertible} elements in $\mathcal{F}$. This is a group for the product $$g(x) u(x)=m_c(g(x), u(x)).$$
Under this law $\mathcal{F}_0$ is a commutative group ($g(x)u(x)=u(x)g(x)$). Again, the identity is the element
$$1 = 1x^0 + 0 x^1 + 0 x^2 + \cdots.$$ The inverse in this group of $g(x)$ is the element $\frac{1}{g(x)}$.

The set $\mathcal{F}_1 = \{\phi(x) \in \mathcal{F} | \phi(x)=\phi_1 x + \phi_2 x^2 + \cdots\}$ is the group of \emph{composable} or \emph{reversible} elements. This is a group for the product given by composition
 $$ \phi(x)\cdot \psi(x) = \phi(x) \circ \psi(x).$$
 The composition product $\circ$ is not commutative, for we have
 $$ (\phi \circ \psi)(x)=\phi(\psi(x)) \ne \psi(\phi(x))= (\psi \circ \phi)(x)$$ in general. The identity element for this group is
 $$id(x)=x=x+0 x^2 + 0 x^3 + \cdots,$$ which can be identified with the element
 $$(0,1,0,0,0,\ldots) \in \mathbb{K}^{\mathbb{N}_0}.$$
 This multiplication defines a map
 $$m_o : \mathcal{F}_1 \times \mathcal{F}_1 \longrightarrow \mathcal{F}_1,$$
 $$    (\phi, \psi) \mapsto \phi \circ \psi,$$
 which as we have seen is not commutative.

If the pair $(g(x), \phi(x))$ defines the infinite lower triangular matrix $M_1$ and the pair $(u(x), \psi(x))$ defines the matrix $M_2$, then the matrix product $M_1 M_2=M_1 \cdot M_2$ is defined by the pair $(f, \theta)$ where
$$ (f, \theta) = (g, \phi) \cdot (u, \psi) = (g.(u\circ \phi), \psi\circ \phi).$$ As for all matrices, this product is associative but not commutative. The identity for this multiplication is the element $(1,x)$. We have
$$ [x^n] x^k = [x^{n-k}] 1 = \delta_{n,k}$$ and hence the matrix corresponding to $(1,x)$ is the (infinite) identity matrix.

\begin{example}
The Riordan array $\left(\frac{1}{1-ax}, \frac{x}{1-ax}\right)$ is represented by the matrix with general element $\binom{n}{k}a^{n-k}$. We have
\begin{eqnarray*}
\left(\frac{1}{1-ax}, \frac{x}{1-ax}\right)\cdot\left(\frac{1}{1-bx}, \frac{x}{1-bx}\right)&=&
\left(\frac{1}{1-ax}\frac{1}{1-b\frac{x}{1-ax}}, \frac{\frac{x}{1-ax}}{1-b\frac{x}{1-ax}}\right)\\
&=&\left(\frac{1}{1-ax} \frac{1-ax}{1-(a+b)x}, \frac{x}{1-(a+b)x}\right)\\
&=& \left(\frac{1}{1-(a+b)x}, \frac{x}{1-(a+b)x}\right).\end{eqnarray*}

This example shows that the set of elements $\left(\frac{1}{1-tx}, \frac{x}{1-tx}\right)$ describes a one-parameter semi-group (in fact, a subgroup) in $\mathcal{R}$. Its infinitesimal generator is given by
$$\left(\begin{array}{ccccccc}
 0 & 0 & 0 & 0 & 0 & 0 & \cdots \\
 1 & 0 & 0 & 0 & 0 & 0 & \cdots \\
 0 & 2 & 0 & 0 & 0 & 0 & \cdots \\
 0 & 0 & 3 & 0 & 0 & 0 & \cdots \\
 0 & 0 & 0 & 4 & 0 & 0 & \cdots \\
 0 & 0 & 0 & 0 & 5 & 0 & \cdots \\
 \vdots & \vdots & \vdots & \vdots & \vdots & \vdots & \ddots \\
\end{array}\right).$$
\end{example}
The inverse matrix $M^{-1}$ of $M$ where $M$ is defined by $(g(x), \phi(x))$ is defined by
$$(g(x), \phi(x))^{-1} = \left(\frac{1}{ (g\circ \bar{\phi})(x)}, \bar{\phi}\right),$$ where
$\bar{\phi}$ is the \emph{reversion} of $\phi(x)$. We sometimes use the notation $\bar{\phi}(x)=\text{Rev}\{\phi\}(x)$.  The existence of $\bar{\phi}$ is guaranteed because
$\phi \in \mathcal{F}_1$. It is the solution $u$ to the equation $\phi(u)=x$ that satisfies $u(0)=0$. Note that
$\bar{x}=x$. Note also that $\overline{\bar{\phi}}(x) = \phi(x)$.

\begin{example} We consider the element $\phi(x)=x(1-x) \in \mathcal{F}_1$. This corresponds to the sequence
$(0,1,-1,0,0,0,\ldots)$. To find $\bar{\phi}(x)$, we solve the equation
$$u(1-u)=x$$ or
$$u^2-u+x=0$$ to get
$$u(x) = \frac{1-\sqrt{1-4x}}{2} \quad\quad \text{or}\quad u(x)=\frac{1+\sqrt{1-4x}}{2}.$$
We obtain $$\bar{\phi}(x)=\frac{1-\sqrt{1-4x}}{2}$$ since we require that the solution satisfy $u(0)=0$.
This is the generating function of the sequence
$$(0,1,1,2,5,14,42,132,\ldots)$$ of the Catalan numbers (with a $0$ pre-pended).
It is easy to show likewise that
$$\text{Rev}\left\{\frac{1-\sqrt{1-4x}}{2}\right\}(x)=x(1-x).$$

The element $(1-x, x(1-x)) \in \mathcal{R}$ is represented by the matrix that has general element given by $(-1)^{n-k}\binom{k+1}{n-k}$, since we have
\begin{eqnarray*}
[x^n](1-x)(x(1-x))^k&=& [x^n]x^k (1-x)^{k+1}\\
&=& [x^{n-k}] \sum_{i=0}^{k+1} \binom{k+1}{i}(-1)^i x^i\\
&=& \binom{k+1}{n-k}(-1)^{n-k}.\end{eqnarray*}

To find the inverse of this element, we introduce the notation
$$c(x) = \frac{1-\sqrt{1-4x}}{2x}.$$
From above we know that
$$ \bar{\phi}(x)=\text{Rev}\{x(1-x)\}(x) = xc(x).$$ The first element of the inverse we seek is given by
$$\frac{1}{(g \circ \bar{\phi})(x)} = \frac{1}{1-xc(x)} = c(x),$$ where we note that the last equality represents a special property of the Catalan numbers.

Thus we get

$$ (1-x, x(1-x))^{-1} = (c(x), xc(x)).$$
\end{example}

This last example provides an example of elements of the \emph{Bell} subgroup $\mathcal{B}$ of the Riordan group $\mathcal{R}$. We have
$$\mathcal{B}= \{ (g, \phi) \in \mathcal{R} | \phi(x)=xg(x) \}.$$

In order to work out a formula for the general element of the matrix corresponding to $(c(x), xc(x))$ we use the Lagrange inversion, in the form given by the Lagrange-B\"urmann theorem. This states that

$$[x^n] G(\bar{\phi}) = \frac{1}{n} [x^{n-1}] G'(x) \left(\frac{x}{\phi}\right)^n, $$
where $G(x) \in \mathbb{K}[[x]]$.

Thus we have

\begin{eqnarray*}
[x^n] c(x)(xc(x))^k &=& [x^n]\frac{1}{x} xc(x)(xc(x))^k\\
&=& [x^{n}] \frac{1}{x} (xc(x))^{k+1} \\
&=& [x^{n+1}] \left(\text{Rev}\{x(1-x)\}(x)\right)^{k+1}\\
&=& \frac{1}{n+1} [x^{n}] (k+1)x^k \left(\frac{x}{x(1-x)}\right)^{n+1}\\
&=& \frac{k+1}{n+1} [x^{n-k}] (1-x)^{-(n+1)}\\
&=& \frac{k+1}{n+1} [x^{n-k}] \sum_{i=0}^{\infty} \binom{-(n+1)}{i}(-1)^i x^i\\
&=& \frac{k+1}{n+1} [x^{n-k}] \sum_{i=0}^{\infty} \binom{n+1+i-1}{i}x^i\\
&=& \frac{k+1}{n+1} \binom{n+n-k}{n-k}\\
&=& \frac{k+1}{n+1} \binom{2n-k}{n-k}.\end{eqnarray*}

We conclude that the matrix corresponding to $(1-x, x(1-x))^{-1}=(c(x),xc(x))$ has its general term given by
$$ \frac{k+1}{n-k+1} \binom{2n-k}{n-k}.$$

 Algebraically, we have

 $$\mathcal{R} = \mathcal{F}_0 \rtimes \mathcal{F}_1, $$ where the symbol $\rtimes$ denotes the semi-direct product.

 This means that we have an exact sequence
 
 \[
 1 \longrightarrow \mathcal{F}_0 \overset{\alpha}{\longrightarrow} \mathcal{R} \overset{\beta}\longrightarrow \mathcal{F}_1 \longrightarrow 1, \]

 where we have the maps
 $$ \alpha : \mathcal{F}_0 \rightarrow \mathcal{R},$$
 $$ g \mapsto (g, x),$$ and
 $$ \beta : \mathcal{R} \rightarrow \mathcal{F}_1,$$
 $$ (g, \phi) \mapsto \phi.$$

 The Riordan group $\mathcal{R}$ has a number of well-studied subgroups. We have already seen the Bell subgroup. The so-called \emph{Appell} subgroup is the subgroup
 $$ \mathfrak{A}=\{ (g(x), \phi(x)) \in \mathcal{R} | \phi(x)=x\},$$
 while the \emph{associated} subgroup or the \emph{Lagrange} subgroup is defined by
 $$\mathfrak{L} = \{ (g(x), \phi(x) \in \mathcal{R} | g(x)=1 \}.$$
 Then $\mathfrak{A}$ is a easily shown to be a normal subgroup of $\mathfrak{R}$.
 Since we have
 $$ ( g(x), x) \cdot (1, \phi(x)) = (g(x).1, (\phi \circ x)(x))=(g(x), \phi(x)) $$ and
$$ (1, \phi(x)) \cdot (g(x), x) = (1.g(x), (x \circ \phi)(x)) = (g(x), \phi(x)) $$  it follows that
$$ \mathcal{R} = \mathfrak{A} \rtimes \mathfrak{L}.$$

In fact, it is clear that we have $$  \mathfrak{A} \cong \mathcal{F}_0 $$ by the mapping
$$ (g(x), x) \mapsto g(x),$$ and we have
$$ \mathfrak{L} \cong \mathcal{F}_1 $$ by the mapping
$$ (1, \phi(x)) \mapsto \phi(x).$$
Thus we re-find that
$$\mathcal{R} \cong \mathfrak{A} \rtimes \mathfrak{L} \cong \mathcal{F}_0 \rtimes \mathcal{F}_1.$$

\section{Tensor products of $R$ modules}

In this section, we let $R$ be a commutative ring with unit,  we let $M_1, M_2, \ldots, M_n$ be a collection of $R$ modules, and we let $A$ be an $R$-module.
Our goal is to define and briefly study the tensor products $M_1 \otimes M_2$,..., $M_1 \otimes \cdots \otimes M_n$. We follow the development in \cite{UnderwoodB} for this.

We shall say that a function $f: M_1 \times M_2 \times \cdots \times M_n \longrightarrow A$ is $R\,n$-linear if for all $i$, $1 \le i \le n$, and all $a_i, a'_i \in M_i, r \in R$ we have
\begin{itemize}
\item $f(a_1, a_2, \ldots, a_i+a'_i, \ldots, a_n)=f(a_1, a_2, \ldots,a_i, \ldots, a_n)+f(a_1, a_2, \ldots, a'_i, \ldots, a_n)$ \\
\item $f(a_1, a_2, \ldots, ra_i,\ldots, a_n) = rf(a_1, a_2, \ldots, a_i, \ldots, a_n)$.
\end{itemize}
For example, an $R$-bilinear map is an $R\,2$-linear map.

A \emph{tensor product} of $M_1, M_2, \ldots, M_n$ over $R$ is an $R$ module $M_1 \otimes M_2 \otimes \cdots \otimes M_n$ together with an $R\,n$-linear map
$$f: M_1 \times M_2 \times \cdots \times M_n \longrightarrow M_1 \otimes M_2 \otimes \cdots \otimes M_n$$ so that for every $R$ module $A$ and $R\,n$-linear map $h: M_1 \times M_2 \times \cdots \times M_n \longrightarrow A$ there exits a unique $R$ module map $\tilde{h}: M_1 \otimes M_2 \otimes \cdots \otimes M_n \longrightarrow A$ for which $\tilde{h}f=h$, that is, the following diagram commutes.

\[
\begin{tikzcd}
M_1 \times M_2 \times \cdots \times M_n \arrow{dr} {f} \arrow{rr}{h} && A \\
& M_1 \otimes M_2 \otimes \cdots \otimes M_n \arrow{ur} {\tilde{h}}
\end{tikzcd}
\]

The tensor product can be constructed as follows. We let $F\langle M_1 \times M_2 \times \cdots \times M_n\rangle $ denote the free $R$-module on the set $M_1 \times M_2 \times \cdots \times M_n$. Let $J$ be the submodule of this $R$-module generated by quantities of the form
$$(a_1, a_2, \ldots, a_i+a'_i, \ldots, a_n)-(a_1,a_2, \ldots, a_i,\ldots, a_n)-(a_1,a_2, \ldots, a'_i,\ldots, a_n),$$
$$(a_1, a_2, \ldots, ra_i,\ldots, a_n)-r(a_1,a_2, \ldots, a_i, \ldots, a_n),$$
for all $i, 1\le i\le n$, and all $a_i,a'_i \in M_i, r \in R$. Let
$$\iota: M_1 \times M_2 \times \cdots \times M_n \longrightarrow F \langle M_1 \times M_2 \times \cdots \times M_n\rangle $$
be the natural inclusion map and let
$$s:F \langle M_1 \times M_2 \times \cdots \times M_n \rangle \longrightarrow F \langle M_1 \times M_2 \times \cdots \times M_n \rangle /J$$ be the canonical surjection. Set $f=s\iota$. Then the quotient space $$F \langle M_1 \times M_2 \times \cdots \times M_n \rangle/J$$ together with the map $f$ is a tensor product which solves the universal mapping problem above \cite{UnderwoodB}.

Furthermore, we have
$$ M_1 \otimes (M_1 \otimes M_2) \cong (M_1 \otimes M_2) \otimes M_3)$$ so that we can use the notation
$$ M_1 \otimes M_2 \otimes M_3 $$ unambiguously. This extends by induction to $M_1 \otimes M_2 \otimes \cdots \otimes M_n$.

Given maps $f_i: M_i \longrightarrow M'_i$ of $R$-modules, for $1 \le i \le n$, there exists a unique map of $R$-modules
$$(f_1\otimes f_2 \cdots \otimes f_n): M_1 \otimes M_2 \otimes \cdots \otimes M_n \rightarrow M'_1 \otimes M'_2 \otimes \cdots \otimes M'_n$$ defined as
$$(f_1\otimes f_2 \cdots \otimes f_n)(a_1\otimes a_2 \otimes \cdots \otimes a_n)=f_1(a_1)\otimes f_2(a_2) \otimes \cdots \otimes f_n(a_n)$$
for all $a_i\in M_i$.

We note the following. If $V_i, 1\le i \le n$ is a finite set of vector spaces over the field $\mathbb{K}$, then $$V_1^* \otimes V_2^* \otimes \cdots \otimes V_n^* \subseteq (V_1 \otimes V_2 \otimes \cdots \otimes V_n)^*,$$ with equality only if all the $V_i$ are finite dimensional.

\section{Algebras, coalgebras and bialgebras}

We begin this section by re-visiting the notion of an algebra, for which we will find useful the notion of the tensor product. Thus we can define a $\mathbb{K}$-algebra to be a triple $(A, m_A, \lambda_A)$ consisting of a vector space $A$ over $\mathbb{K}$, and $\mathbb{K}$-linear maps $m_A: A \otimes A \rightarrow A$ and $\lambda_A: \mathbb{K} \rightarrow A$ that satisfy the following conditions.
\begin{itemize}
\item We have a commutative diagram:
\[
\begin{tikzcd}
A \otimes A \otimes A \arrow{d} {m_A \otimes I_A} \arrow{r} {I_A \otimes m_A} & A \otimes A \arrow{d} {m_A} \\
A \otimes A \arrow{r} {m_A} & A
\end{tikzcd}
\]
Here, $I_A$ is the identity on $A$, and we have $I_A \otimes m_A: A \otimes A \otimes A \rightarrow A \otimes A$ defined by $a\otimes b \otimes c \mapsto a \otimes m_A(b \otimes c)$ and similarly for $m_A \otimes I_A$. The commutativity of the diagram thus means that for all $a, b, c \in A$, we have
\begin{equation} m_A(I_A \otimes m_A)(a \otimes b \otimes c) = m_A(m_A \otimes I_A)(a \otimes b \otimes c).\label{ass}\end{equation}
\item We have a commutative diagram
\[
\begin{tikzcd}
A \otimes  \mathbb{K} \arrow[d, "s_2"'] \arrow{r} {I_A \otimes  \lambda_A} & A \otimes A  \arrow{dl} {m_A} \\
A  & \mathbb{K} \otimes A  \arrow[u, "\lambda_A \otimes I_A"'] \arrow{l} {s_1}
\end{tikzcd}
\]
Here, the map $s_1: \mathbb{K} \otimes A \rightarrow A$ is defined by $r \otimes a \mapsto ra$ and the map $s_2: A \otimes \mathbb{K} \rightarrow A$ is defined by $a \otimes r \mapsto ra$.  We have
\begin{equation}m_A(I_A \otimes \lambda_A)(a \otimes r)=ra=m_A(\lambda_A \otimes I_A)(r \otimes a).\label{unit}\end{equation}
\end{itemize}

The map $m_A$ is the multiplication map (for $A$) and the map $\lambda_A$ is the \emph{unit map}. The maps $s_1$ and $s_2$ represent scalar multiplication. The property (\ref{ass}) is called the \emph{associative property} and the property (\ref{unit}) is called the \emph{unit property}.

It is straightforward to show that this definition of a $\mathbb{K}$-algebra coincides with that given before.

\begin{example} The polynomial $\mathbb{K}[x]$ is a $\mathbb{K}$-algebra with multiplication
$$m_{\mathbb{K}} : \mathbb{K}[x] \otimes \mathbb{K}[x] \rightarrow \mathbb{K}[x]$$ given by ordinary polynomial multiplication, and $\lambda_{\mathbb{K}[x]} : \mathbb{K} \rightarrow \mathbb{K}[x]$ defined as
$r \mapsto r.1$, for all $r \in \mathbb{K}$.
\end{example}

\begin{example} The ring of power series $\mathbb{K}[[x]]$ over $\mathbb{K}$ is a $\mathbb{K}$-algebra. The multiplication is given by
  $$m_{\mathbb{K}} : \mathbb{K}[[x]] \otimes \mathbb{K}[[x]] \rightarrow \mathbb{K}[[x]]$$
  $$ f \otimes g \mapsto \sum_{n=0}^{\infty} (\sum_{k=0}^n f_k g_{n-k}) x^n.$$
The unit $\lambda_{\mathbb{K}[[x]]} : \mathbb{K} \rightarrow \mathbb{K}[[x]]$ is defined as
$r \mapsto r.1$, for all $r \in \mathbb{K}$.
\end{example}

These algebras are commutative, where we define a $\mathbb{K}$-algebra  $A$ to be commutative if we have
$$ m_A \tau = m_A $$ where $\tau$ denotes the \emph{twist map} defined by $\tau(a \otimes b)=b \otimes a$ for $a, b \in A$.

An example of a non-commutative $\mathbb{K}$-algebra is the algebra $\mathcal{F}_1$ of composable power series.

Let $(A, m_A, \lambda_A), (B, m_B, \lambda_B)$ be two $\mathbb{K}$-algebras. A $\mathbb{K}$-\emph{algebra homomorphism} from $A$ to $B$ is a map of additive groups $\phi: A \rightarrow B$ for which $\phi(1_A)=1_B$  with $$ \phi(m_A(a\otimes a'))=m_B(\phi(a) \otimes \phi(a')),$$ and
$$ \phi(\lambda_A(r))=\lambda_B(r)$$ for $a,a' \in A, r \in \mathbb{K}$.

The usefulness of using commutative diagrams to describe $\mathbb{K}$ will now come into play when we define the notion of \emph{coalgebra}. Essentially we will need to reverse the arrows in some of our diagrams.

 Thus a $\mathbb{K}$-coalgebra is a triple $(C, \Delta_C, \epsilon_C)$ consisting of a vector space $C$ over $\mathbb{K}$ and $\mathbb{K}$-linear maps
$$\Delta_C : C \longrightarrow C \otimes C$$ and
$$ \epsilon_C : C \longrightarrow \mathbb{K}$$ that satisfy the following conditions.
\begin{itemize}
\item The following diagram commutes.
\[
\begin{tikzcd}
C \arrow{r} {\Delta_C} \arrow{d} {\Delta_C} & C \otimes C \arrow{d} {I_C \otimes \Delta_C}\\
C \otimes C \arrow{r} {\Delta_C \otimes I_C} & C \otimes C \otimes C
\end{tikzcd}
\]
Here the map $I_C:C \rightarrow C$ is the identity map and the maps $I_C \otimes \Delta_C: C\otimes C \rightarrow C\otimes C \otimes C$ and $\Delta_C\otimes I_C:C \otimes \rightarrow C \otimes C\otimes C$ are defined by $a \otimes b \mapsto a \otimes \Delta_C(b)$ and $a \otimes b \mapsto \Delta_C(a)\otimes b$, for all $a, b \in C$, respectively.
Thus for all $c \in C$, we have
\begin{equation} (I_C \otimes \Delta_C) \Delta_C (c) = (\Delta_C \otimes I_C) \Delta_C(c).\label{coass} \end{equation}
 \item The following diagram commutes.
 \[
 \begin{tikzcd}
 C \arrow{d} {- \otimes 1} \arrow{rd} {\Delta_C} \arrow{r} {1 \otimes -} & \mathbb{K}\otimes C \\
 C \otimes \mathbb{K} & C \otimes C \arrow{l} {I_C \otimes \epsilon_C} \arrow {u} {\epsilon_C \otimes I_C}
 \end{tikzcd}
 \]
 \end{itemize}

 Here the maps $- \otimes 1$ and $1 \otimes -$ are defined by $c \mapsto c \otimes 1$ and $c \mapsto 1 \otimes c$, respectively. Equivalently, we have
\begin{equation} (\epsilon_C \otimes I_C)\Delta_C = 1 \otimes c, \quad\quad (I_C \otimes \epsilon_C) \Delta_C(c)=c \otimes 1. \label{counit} \end{equation}

The map $\Delta_C$ is called the \emph{comultiplication} map. The map $\epsilon_C$ is called the \emph{counit} map.  The condition (\ref{coass}) is called the \emph{coassociative property} and the condition (\ref{counit}) is called the \emph{counit property}.

A $\mathbb{K}$-coalgebra $C$ is \emph{cocommutative} if $$\tau(\Delta_C)(c))=\Delta_C(c)$$ for all $c \in C$.

\begin{example} Let $\mathbb{K}[x]$  denote the $\mathbb{K}$-vector space of polynomials in the indeterminate $x$. Let
$$ \Delta_{\mathbb{K}[x]} : \mathbb{K}[x] \longrightarrow \mathbb{K}[x] \otimes \mathbb{K}[x]$$ be the $\mathbb{K}$-linear map defined on the $\mathbb{K}$-basis $\{1,x, x^2,\ldots\}$ as
$$ \Delta_{\mathbb{K}[x]}(x^m) = x^m \otimes x^m,$$ and
let
$$ \epsilon_{\mathbb{K}[x]} : \mathbb{K}[x] \longrightarrow \mathbb{K}$$ be the $\mathbb{K}$-linear map defined on $\{1,x, x^2,\ldots\}$ as
$$ \epsilon_{\mathbb{K}[x]}(x^n) = 1.$$
Then the triple $(\mathbb{K}[x], \Delta_{\mathbb{K}[x]}, \epsilon_{\mathbb{K}[x]})$ is a $\mathbb{K}$-coalgebra. 
\end{example}

\begin{example} We follow \cite{BFK} in this example, which is built around the group algebra or coordinate ring of the set $\mathcal{F}_0$. Thus we denote by $\mathbb{C}(\mathcal{F}_0)$ the set of functions
$$ F: \mathcal{F}_0 \longrightarrow \mathbb{C}$$ that are polynomial with respect to an appropriate basis. Choosing the basis $\{1=[x^0], [x], [x^2], \ldots\}$ we obtain the identification
$$ \mathbb{C}(\mathcal{F}_0) \cong \mathbb{C}[[x], [x^2], [x^3], \ldots].$$
We then obtain a map
$$ \Delta_0 : \mathbb{C}(\mathcal{F}_0) \longrightarrow \mathbb{C}(\mathcal{F}_0) \otimes \mathbb{C}(\mathcal{F}_0),$$
$$ [x^n] \mapsto (f \otimes g \mapsto [x^n]fg).$$
Since we have $[x^n]fg = \sum_{k=0}^n [x^k] f [x^{n-k}]g$, we obtain that
$$ \Delta_0 [x^n] = \sum_{k=0}^n [x^k] \otimes [x^{n-k}].$$

We define a counit $\epsilon_0$ on $\mathbb{C}(\mathcal{F}_0)$ by
$$ \epsilon_0 : \mathbb{C}(\mathcal{F}_0) \longrightarrow \mathbb{C},$$
$$               [x^n] \mapsto [x^n]1 $$
Thus we have $\epsilon_0([x^n]) = \delta_{n,0}$.
Then $(\mathbb{C}(\mathcal{F}_0), \Delta_0, \epsilon_0)$ is a coalgebra. As a ring of polynomials, with multiplication of polynomials as the multiplication, and the unit element $[x^0]=1$, the set $\mathbb{C}(\mathcal{F}_0)$ is of course an algebra as well.
\end{example}

\begin{example} We again follow \cite{BFK} in this example, which is built around the coordinate ring of the set $\mathcal{F}_1$. For this, we define the family of functionals $a_n$ on $\mathcal{F}_0$ by
$$a_n(\phi) = \frac{1}{(n+1)!}\frac{d^{n+1}}{dx^{n+1}}\phi|_{x=0}.$$ In other words, $a_n=\frac{1}{(n+1)!}[x^{n+1}]$. We identify the coordinate ring $\mathbb{C}(\mathcal{F}_1)$ with the polynomial ring $\mathbb{C}[a_1, a_2, \ldots]$ in infinitely many variables $a_1, a_2, \ldots$. We can define a co-product on $\mathbb{C}(\mathcal{F}_1)$ as follows:
$$ \Delta_1 : \mathbb{C}(\mathcal{F}_1) \longrightarrow \mathbb{C}(\mathcal{F}_1) \otimes \mathbb{C}(\mathcal{F}_1),$$
$$ a_n \mapsto (\phi \otimes \psi \mapsto a_n(\phi \circ \psi)).$$ Thus the co-product for the generators of
$\mathbb{C}(\mathcal{F}_1)$ can be extracted from the standard duality condition
$$\langle \Delta_1 a_n, \phi \otimes \psi\rangle = a_n(\phi \circ \psi),$$
where $\langle a_n, \phi \rangle = a_n(\phi)$ and $\langle a_n \otimes a_m, \phi \otimes \psi\rangle = a_n(\phi)a_m(\psi)$.
\end{example}

Let $C$  be a $\mathbb{K}$-coalgebra. A non-zero element $c$ of $C$ for which $\Delta_C(c)=c \otimes c$ is called a \emph{grouplike} element of $C$. Necessarily we then have $\epsilon_C(c)=1$. The set of grouplike elements $G(C)$ of $C$ is a linearly independent subset of $C$.

Now let $C, D$ be coalgebras over $\mathbb{K}$. A $\mathbb{K}$-linear map $\phi: C \rightarrow D$ is a \emph{coalgebra homomorphism} if
$$ (\phi \otimes \phi)\Delta_C(c)=\Delta_D(\phi(c))$$ and
$$ \epsilon_C(c) = \epsilon_D(\phi(c))$$ for all $c \in C$.

\begin{example} The field $\mathbb{K}$ as a vector space over itself is a $\mathbb{K}$-coalgebra where the comumltiplicaton map $\Delta_{\mathbb{K}}: \mathbb{K} \rightarrow \mathbb{K} \otimes \mathbb{K}$ is defined by
$\Delta_{\mathbb{K}}(a) = a \otimes 1$ and the counit map $\epsilon_{\mathbb{K}}: \mathbb{K}\rightarrow \mathbb{K}$ is defined by $\epsilon_{\mathbb{K}}(a)=a$. This is the \emph{trivial coalgebra}.

If now $C$ is a $\mathbb{K}$-coalgebra, then the counit map $\epsilon_C: C \rightarrow \mathbb{K}$ is a homomorphism of $\mathbb{K}$ algebras.
\end{example}

If $\phi:  C \rightarrow D$ is a homomorphism of coalgebras, and if $c$ is a grouplike element of $C$, then $\phi(c)$ is a grouplike element of $D$.

A coalgebra homomorphism $\phi: C \rightarrow D$ that is injective and surjective is an \emph{isomorphism of coalgebras}.

Duality will play an important role in the sequel. Thus we look at this in the context of algebras and coalgebras. If $C$ is a $\mathbb{K}$-coalgebra, then we denote by $C^*$ its linear dual. The important fact now is that if $C$ is a $\mathbb{K}$-coalgebra, then its dual $C^*$ is an algebra. Under this duality, the dual mapping
$$ \Delta_C^* : (C \otimes C)^* \longrightarrow C^*$$ restricts to a $\mathbb{K}$-linear map $m_{C^*}$ to $C^* \otimes C^* \subseteq (C \otimes C)^*$ defined as
\begin{eqnarray*} m_{C^*}(f \otimes g)(c)&=&\Delta_C^*(f \otimes g)(c)\\
&=& (f \otimes g) (\Delta_C(c))\\
&=&= \sum_{(c)} f(c_{(1)})g(c_{(2)}).\end{eqnarray*}

The transpose of the counit map of $C$ is
$$\epsilon_C^* : \mathbb{K}^*=\mathbb{K} \rightarrow C^*$$ defined by
$$\epsilon_C^*(r)(c)=r(\epsilon_C(c))=r\epsilon_C(c)$$ for $r \in \mathbb{K}, c\in C$.  We set
$\lambda_{C^*}=\epsilon_C^*$ and define maps $$I_{C^*}\otimes \lambda_{C^*}: C^* \otimes \mathbb{K} \rightarrow C^* \otimes C^*,$$
$$ f \otimes r \mapsto f \otimes \lambda_{C^*}(r),$$ and
$$\lambda_{C^*} \otimes I_{C^*} : \mathbb{K} \otimes C^* \rightarrow C^* \otimes C^*,$$
$$ r \otimes f \mapsto \lambda_{C^*}(r) \otimes f,$$ for $f\in C^*, r \in \mathbb{K}$.

With these definitions, it can be shown that if $(C, \Delta_C, \epsilon_C)$ is a coalgebra, then $(C^*, m_{C^*}, \lambda_{C^*})$ is an algebra.

We may ask if the converse is also true, that is, if $(A, m_A, \lambda_A)$ is an algebra, does $A^*$ have the structure of a coalgebra? The transpose of the multiplication map $m_A^*$ is such that $$m_A^* : A^* \rightarrow (A \otimes A)^*,$$ but in the infinite dimensional case, we have that $A^* \otimes A^*$ is a proper subset of $(A \otimes A)^*$, and so we may not have a transpose mapping $A^* \rightarrow A^* \otimes A^*$ which would be a necessary condition for $A^*$ to be an algebra. To overcome this difficulty, we can proceed as follows. For $A$ a $\mathbb{K}$-algebra, we define the finite dual $A^o$ of $A$ by
$$ A^o = \{ f \in A^* | f\,\text{ vanishes on some ideal}\, I \subseteq A\,\text{ of finite codimension}\}.$$
If $A$ is an algebra, then it can be shown that $A^o$ is a coalgebra.
\begin{example} The collection of $k$-th order linearly recursive sequences over $\mathbb{K}$ of all orders $k>0$ can be identified with the finite dual $\mathbb{K}[x]^o$.
\end{example}

Sometimes, the two structures, algebra and coalgebra, can co-exist. A $\mathbb{K}$-bialgebra is a $\mathbb{K}$ vector space $B$ together with maps $m_B, \lambda_B, \Delta_B, \epsilon_B$ that satisfy the following conditions:
\begin{enumerate}
\item $(B, m_B, \lambda_B)$ is a $\mathbb{K}$-algebra and $(B, \Delta_B, \epsilon_B)$ is a $\mathbb{K}$-algebra, \\
\item $\Delta_B$ and $\epsilon_B$ are homomorphisms of $\mathbb{K}$-algebras.
\end{enumerate}

If $B$ is a bialgebra, it can be shown that $B^o$ will also be a bialgebra.

\begin{example} The sets $\mathbb{C}(\mathcal{F}_0)$ and $\mathbb{C}(\mathcal{F}_1)$ are bialgebras.
\end{example}

\section{Hopf algebras and the Riordan group}

We are now in a position to define what we mean by a Hopf algebra. Thus a $\mathbb{K}$-Hopf algebra is a bialgebra over a field $\mathbb{K}$
$$H = \{H, m_H, \lambda_H, \Delta_H, \epsilon_H\}$$ together with a $\mathbb{K}$-linear map
$$ \sigma_H : H \longrightarrow H$$
that satisfies
\begin{equation} m_H(I_H \otimes \sigma_H) \Delta_H(h) = \epsilon_H(h)1_H = m_H(\sigma_H \otimes I_H)\Delta_H(h) \label{coin}\end{equation} for all $h \in H$.

The map $\sigma_H$ is called the \emph{antipode} or \emph{coinverse} and property (\ref{coin}) is called the antipode or the coinverse property.

\begin{example} Consider the coordinate ring $\mathbb{C}(\mathcal{F}_0)$ of the set of invertible power series $\mathcal{F}_0$ over $\mathbb{C}$. We can define a coinverse map
$$S_0 : \mathbb{C}(\mathcal{F}_0) \longrightarrow \mathbb{C}(\mathcal{F}_0),$$
$$[x^n] \mapsto (f \mapsto [x^n] \frac{1}{f}).$$
To show that this is a coinverse map, we must show that
$$\sum_{k=0}^n [x^k]  \otimes (f\mapsto [x^{n-k}]\frac{1}{f}) = \sum_{k=0}^n (f \mapsto[x^k]\frac{1}{f})\otimes[x^{n-k}].$$
But this is true since
$$ \delta_{n,0}=[x^n]1 = [x^n] f.\frac{1}{f} = \sum_{k=0}^n [x^k]f[x^{n-k}]\frac{1}{f}$$ and
$$ \delta_{n,0}=[x^n]1 = [x^n] \frac{1}{f}.f = \sum_{k=0}^n [x^k]\frac{1}{f} [x^{n-k}] f.$$
We have already seen that $\mathbb{C}(\mathcal{F}_0)$ is a bialgebra. With this coinverse or antipode mapping, it can be shown that $\mathcal{H}_0=\mathbb{C}(\mathcal{F}_0)$ is a Hopf algebra \cite{BFK}.

As shown in \cite{BFK}, it is possible to recover $\mathcal{F}_0$ from $\mathcal{H}_0$.
Thus we have that
$$ \mathcal{F}_0 \cong Hom_{Alg}(\mathcal{H}_0, \mathbb{C}),$$ where this last expression denotes the group of algebra homomorphisms (or characters) on  $\mathcal{H}_0$, with the convolution product defined on the generators by
$$ (\alpha \beta)[x^n]:= m \circ (\alpha \otimes \beta) \circ \Delta_0 [x^n],$$ for any algebra homomorphisms $\alpha, \beta$ on $\mathcal{H}_0$. Here, $m$ is multiplication on $\mathbb{C}$. 

\[
\begin{tikzcd}
\mathcal{H}_0  \arrow{r} {\Delta_0} \arrow{ddr}{\alpha \beta} & \mathcal{H}_0 \otimes \mathcal{H}_0 \arrow{d} {\alpha \otimes \beta} \\
& \mathbb{C} \otimes \mathbb{C} \arrow{d} {m} \\
& \mathbb{C}
\end{tikzcd}
\]

The isomorphism is given by
$$ \mathcal{F}_0 \cong Hom_{Alg}(\mathcal{H}_0, \mathbb{C}),$$
$$ f \mapsto \alpha_f := ([x^n] \mapsto [x^n]f=f_n).$$
\end{example}

\begin{example} In this example, we consider the coordinate ring $\mathbb{C}(\mathcal{F}_1)$ of the set of composable power series over $\mathbb{C}$. This set is a bialgebra. We can define a coinverse map
$$S_1 : \mathbb{C}(\mathcal{F}_1) \longrightarrow \mathbb{C}(\mathcal{F}_1),$$
$$ a_n \mapsto (f \mapsto a_n(\bar{f})).$$
 With this coinverse or antipode mapping, it can be shown that $\mathcal{H}_1=\mathbb{C}(\mathcal{F}_1)$ is a Hopf algebra \cite{BFK}.

Again, it is possible to recover $\mathcal{F}_1$ from $\mathcal{H}_1$.
Thus we have that
$$ \mathcal{F}_1 \cong Hom_{Alg}(\mathcal{H}_1, \mathbb{C}),$$ where this last expression denotes the group of algebra homomorphisms (or characters) on  $\mathcal{H}_0$, with the convolution product defined on the generators by
$$ (\alpha \beta) a_n := m \circ (\alpha \otimes \beta) \circ \Delta_1 a_n,$$ for any algebra homomorphisms $\alpha, \beta$ on $\mathcal{H}_0$. Here, $m$ is multiplication on $\mathbb{C}$.

\[
\begin{tikzcd}
\mathcal{H}_1  \arrow{r} {\Delta_1} \arrow{ddr}{\alpha \beta} & \mathcal{H}_1 \otimes \mathcal{H}_1 \arrow{d} {\alpha \otimes \beta} \\
& \mathbb{C} \otimes \mathbb{C} \arrow{d} {m} \\
& \mathbb{C}
\end{tikzcd}
\]

 The isomorphism is given by
$$ \mathcal{F}_1 \cong Hom_{Alg}(\mathcal{H}_0, \mathbb{C}),$$
$$ \phi \mapsto \alpha_{\phi} := (a_n \mapsto a_n(\phi)).$$
\end{example}

We now observe that we have an isomorphism of algebras

$$ \mathbb{C}(\mathcal{F}_0 \ltimes \mathcal{F}_1) \cong \mathcal{H}_0 \otimes \mathcal{H}_1,$$ or
$$ \mathbb{C}(\mathcal{R}) \cong \mathcal{H}_0 \otimes \mathcal{H}_1.$$

It is clear that the relationship between $\mathbb{C}(\mathcal{R})$ and the semi-direct product $\mathcal{H}_0 \rtimes \mathcal{H}_1$ \cite{Molnar} deserves further research.

\section{Conclusions}

To date, the Riordan group has been applied, in the main, for investigations in the area of combinatorics. Its appearance in other areas, all associated to its Lie group nature, and in particular where the emphasis is in applications to mathematical physics \cite{BFK, Goodenough}, indicates that it may be important to study the group from different perspectives. This note revolves around the perspective of Hopf algebras, though it barely scratches the surface. 

We have looked only at the so-called ``ordinary'' Riordan group, defined using ordinary generating functions. Of equal importance is the ``exponential'' Riordan group, defined using exponential generating functions \cite{FBK, Goodenough}. The Faa di Bruno formula comes into play in the definition of the relevant co-product in this context. The corresponding power series are often referred to as ``divided power'' series in the mathematical physics area. Depending on the area of application, it would be useful to develop a full theory of Hopf algebra association for the general Riordan group as describe in \cite{Wang}.

The paper \cite{BFK} indicates how the Riordan group might be generalized using invertible series with non-commutative coefficients. This would further motivate studies of the semi-direct product $\mathcal{H}_0 \rtimes \mathcal{H}_1$ and its generalizations.

\section {Appendix}

Given its importance in the discussion, we summarize the rules of operation of the operator $[x^n]$ \cite{MC}.
\begin{center}
\fbox{\begin{minipage}{16 cm}
\begin{tabular}{cllcl}\\
MC1 & Linearity & $[x^n] (r f(x) + s g(x))$ & $=$ & $r [x^n] f(x) + s [x^n] g(x)$\\ \\
MC2 & Shifting & $[x^n]x f(x)$ & $=$ & $ [x^{n-1}] f(x)$ \\ \\
MC3 & Differentiation & $[x^n] f'(x)$ & $=$ & $(n+1)[x^{n+1}] f(x) $ \\ \\
MC4 & Convolution & $ [x^n] g(x)f(x)$ & $=$ & $\sum_{k=0}^n ([x^k]g(x))[x^{n-k}]f(x)$ \\ \\
MC5 & Composition & $ [x^n] g(f(x)) $ & $=$ & $\sum_{k=0}^{\infty} ([x]^k g(x)) [x^n]f(x)^k $\\ \\
MC6 & Inversion & $ [x^n] \bar{f}(x)^k $ & $=$ & $ \frac{k}{n} [x^{n-k}] \left(\frac{x}{f(x)}\right)^n$
\end{tabular}
\end{minipage}}
\end{center}

Note that in (MC1), $r, s \in R$. We can extend rule (MC2) to the following.
\begin{center}
\fbox{\begin{minipage}{5 cm}
\begin{center}
$$ [x^n] x^k f(x) = [x^{n-k}] f(x).$$
\end{center}
\end{minipage}}
\end{center}
There is a more general form of rule (MC6), which is known as \emph{Lagrange Inversion}\index{Lagrange inversion}. We have
\begin{center}
\fbox{\begin{minipage}{6 cm}
\begin{center}
$$[x^n] G(\bar{f})= \frac{1}{n}[x^{n-1}]G'(x)\left(\frac{x}{f}\right)^n,$$
where $G(x) \in R[[x]]$.
\end{center}
\end{minipage}}
\end{center}

\bigskip
\hrule
\bigskip
\noindent 2010 {\it Mathematics Subject Classification}: Primary
33B10; Secondary 33B20, 16T05, 05A15
\noindent \emph{Keywords:} Power series, Riordan group, Hopf algebra

\end{document}